\documentclass[conference,letterpaper]{IEEEtran}
\IEEEoverridecommandlockouts

\usepackage[utf8]{inputenc} 
\usepackage[T1]{fontenc}
\usepackage{url}
\usepackage{ifthen}
\usepackage{cite}
\usepackage[cmex10]{amsmath}
\usepackage{amsmath}
\usepackage{amsthm}
\usepackage{amsfonts}
\usepackage{amssymb}
\usepackage{mathtools}
\usepackage{color}
\usepackage{verbatim}
\usepackage{mathtools}
\usepackage{mathrsfs}
\usepackage[normalem]{ulem}
\usepackage{enumitem}
\usepackage{hyperref}
\usepackage[ruled,vlined]{algorithm2e}
\usepackage{algpseudocode}

\definecolor{purple}{rgb}{0.5, 0.0, 0.5}
\definecolor{dark_green}{rgb}{0.0, 0.5, 0.0}

\newcommand{\cyan}{\color{cyan}}
\newcommand{\purple}{\color{purple}}
\newcommand{\white}{\color{white}}

\newcommand{\ow}{\mathcal{O}}

\newcommand{\I}{\mathcal{I}}

\newcommand{\ub}{\bold{u}}

\newcommand{\vb}{\bold{v}}

\newcommand{\yb}{\bold{y}}

\newcommand{\Gb}{\bold{G}}
\newcommand{\Gbt}{\tilde{\bold{G}}}

\newcommand{\texT}{\mathrm{T}}

\newcommand{\R}{\mathbb{R}}

\newcommand{\N}{\mathbb{N}}

\newcommand{\gb}{\bold{g}}

\newcommand{\Ab}{\bold{A}}
\newcommand{\Acal}{\mathcal{A}}

\newcommand{\bb}{\bold{b}}
\newcommand{\bbh}{\hat{\bold{b}}}

\newcommand{\Bb}{\bold{B}}
\newcommand{\ab}{\bold{a}}
\newcommand{\abt}{\tilde{\bold{a}}}

\newcommand{\cb}{\bold{c}}

\newcommand{\Eb}{\bold{E}}
\newcommand{\eb}{\bold{e}}
\newcommand{\Ib}{\bold{I}}

\DeclarePairedDelimiter\floor{\lfloor}{\rfloor}

\newtheorem{Thm}{Theorem}
\newtheorem{Cor}[Thm]{Corollary}
\newtheorem{Prop}[Thm]{Proposition}

\newcommand{\ind}{\text{\color{white}.$\quad$}}

\begin{document}
\title{Straggler Robust Distributed Matrix Inverse and Pseudoinverse Approximation}

\author{
  \IEEEauthorblockN{$\textbf{Neophytos Charalambides}^{\natural}$, $\textbf{Mert Pilanci}^{\sharp}$, \textbf{and} $\textbf{Alfred O. Hero III}^{\natural}$}
  \thanks{This work was partially supported by grant ARO W911NF-15-1-0479. We also thank Nima Mohseni for a suggestion in Section \ref{pr_pinv_alg_subsec}.}
  \IEEEauthorblockA{$\text{\white.}^{\natural}$EECS Department University of Michigan $\text{\white.}^{\sharp}$EE Department Stanford University\\
  Email: neochara@umich.edu, pilanci@stanford.edu, hero@umich.edu}
\vspace{-4mm}
}

\maketitle

\begin{abstract}
A cumbersome operation in signal processing, numerical analysis and linear algebra, optimization and machine learning, is inverting large full-rank matrices. We propose an approximate matrix inversion algorithm which uses a black-box least squares optimization solver as a subroutine, to give an estimate of the inverse and pseudoinverse of real full-rank matrices. The proposed algorithms can be implemented in a distributed manner by splitting the workload column-by-column, and allocating the computations to multiple workers. Furthermore, it does not require a matrix factorization, e.g. LU, SVD or QR decomposition beforehand. We further incorporate our algorithms in a numerically stable binary repetition coded computing scheme, which makes them robust to \textit{stragglers}, i.e. non-responsive workers.
\end{abstract}

\section{Introduction}
\label{intro}

Inverting a matrix is common operation in numerous applications in domains such as social networks, numerical analysis and integration, machine learning, and scientific computing \cite{GS59,Hig02}. It is one of the most important operations, as it reverses a system. A common way of inverting a matrix is by performing Gaussian elimination, which in general takes $\ow(n^3)$ operations for square matrices of order $n$. In high-dimensional applications, this is cumbersome.

More elegant algorithms with lower complexity exist, which are similar to matrix multiplication algorithms. Some popular algorithms are the ones by Strassen with complexity $\ow(n^{2.807})$ \cite{Str69}, Coppersmith-Winograd with $\ow(n^{2.375477})$ \cite{CW87}, and Le Gall with $\ow(n^{2.3728639})$ \cite{Gal14}. Of these, the Strassen algorithm is used more often. Many other inversion algorithms assume specific structure on the matrix, require a matrix-matrix product, or use a matrix factorization, e.g. LU, SVD or QR decomposition \cite{TB97}. Methods for matrix inversion or factorization are often referred to as \textit{direct methods}, in contrast to \textit{iterative methods}, which gradually converge to the solution \cite{DRSL16,PV20}. The most computationally efficient direct methods compute some form of the inverse, and are asymptotically equivalent. These have complexity $\ow(n^{\omega})$ for $\omega<2.376$.

Due to the prevalence of high-dimensional datasets, distributed algorithms have been of interest, where a network of workers perform certain subtasks in parallel \cite{MHGB18,LLA16,XMA14,QQSV01,QQSV98,YLX13,LKV96,BF88}. Some drawbacks of these algorithms are that they make assumptions on the matrix structure and assume distributed memory, they are specific for distributed and parallel computing platforms (e.g. Apache Spark and Hadoop, MapReduce, CUDA), or require a matrix factorization (e.g. LU, QR). They may also require heavy and multiple communication instances, which makes these algorithms unsuitable for iterative methods.

In this paper, we propose a centralized distributed algorithm which \textit{approximates} the inverse of a nonsingular $\Ab\in\R^{n\times n}$ and makes none of the aforementioned assumptions. The main idea behind the algorithm is that the worker nodes use a least squares solver to approximate one or more columns of $\Ab^{-1}$, which the central server then concatenates to obtain the approximation $\widehat{\Ab^{-1}}$. While other iterative procedures are applicable, we use a steepest descent (SD) for our least squares solver, and also present simulation results with the conjugate gradient method (CG). The proposed algorithm is simple, and can be made robust to \textit{stragglers}; nodes which fail to compute their task or have longer response time than others, through a coded computing scheme. We extend this idea to distributed approximation of the pseudoinverse $\Ab^{\dagger}$ for $\Ab$ full-rank.

The coded computing scheme we propose which guarantees straggler resiliency is an adaptation of the binary gradient code from \cite{CMH20}, a numerically stable repetition scheme with efficient decoding. For our scheme, we set up a framework analogous to that proposed for gradient coding \cite{TLDK17}.

We point out two articles which have similarities to the work presented in this paper \cite{GHSTM11,YGK17}. Firstly, our approach to inverting $\Ab$ is similar in nature to \cite{GHSTM11}, which uses stochastic gradient descent to approximate matrix factorizations distributively. Secondly, the formulation of our underlying optimization problem: minimize $\|\Ab\Bb-\Ib_n\|_F^2$ by estimating the columns of $\Bb$, is equivalent to the problem studied in \cite{YGK17}, which deals with approximating linear inverse problems in the presence of stragglers. The drawbacks of the coded computing scheme provided in \cite{YGK17}, is that it is geared towards specific applications (e.g. personalized PageRank), makes assumptions on the covariance between the signals comprising the linear system and the accuracy of the workers, and assumes an additive decomposition of $\Ab$. Furthermore, the approximation algorithm in \cite{YGK17} is probabilistic. We on the other hand make \textit{no} assumption on $\Ab$ other than the fact that it is non-singular, and our algorithm is not probabilistic.

The paper is organized as follows. In Section \ref{prel_sec} we recall basics of matrix inversion, least squares approximation and SD. In Section \ref{Appr_alg_sec} we propose our matrix inverse and pseudoinverse algorithms. In Section \ref{Str_pr_sec} we present a binary repetition scheme which makes our algorithms robust to stragglers. Finally, in Section \ref{exper_sec} we present some numerical simulations on randomly generated matrices. Our main contributions are:
\begin{itemize}[noitemsep,nolistsep]
  \item A new matrix inverse approximation algorithm and associated pseudoinverse approximation algorithm.
  \item Theoretical guarantee on the approximation errors.
  \item Development of a binary coded computing scheme, which makes the proposed algorithms robust to stragglers.
  \item Experimental justification of our theory.
\end{itemize}

\section{Preliminary Background}
\label{prel_sec}

Recall that a nonsingular matrix is a square matrix $\Ab\in\R^{n\times n}$ of full-rank, which has a unique inverse $\Ab^{-1}$ such that $\Ab\Ab^{-1}=\Ab^{-1}\Ab=\Ib_n$. The simplest way of computing $\Ab^{-1}$ is by performing Gaussian elimination on $\big[\Ab|\Ib_n\big]$, which gives $\big[\Ib_n\big|\Ab^{-1}]$ in $\ow(n^3)$ operations. In Algorithm \ref{inv_alg}, we approximate $\Ab^{-1}$ column-by-column.

For full-rank rectangular matrices $\Ab\in\R^{n\times m}$ where $n>m$, one resorts to the left Moore–Penrose pseudoinverse $\Ab^{\dagger}\in\R^{m\times n}$, for which $\Ab^{\dagger}\Ab=\Ib_m$. In Algorithm \ref{pinv_alg}, we present how to approximate the left pseudoinverse of $\Ab$, by using the fact that $\Ab^{\dagger}=(\Ab^T\Ab)^{-1}\Ab^T$; since $\Ab$ is full-rank. The right pseudoinverse $\Ab^{\dagger}=\Ab^T(\Ab\Ab^T)^{-1}$ of $\Ab^{m\times n}$ where $m<n$, can be obtained by a modification of Algorithm \ref{pinv_alg}.

In the proposed algorithms we approximate $n$ instances of the least squares minimization problem
\begin{equation}
\label{OLS}  
  \theta^{\star}_{ls} = \arg\min_{\theta\in\R^m} \left\{\|\Ab\theta-\yb\|_2^2\right\} 
\end{equation}
for $\Ab\in\R^{n\times m}$ and $\yb\in\R^n$. In many applications $n\gg m$, where the rows represent the feature vectors of a dataset. This has the closed-form solution $\theta^{\star}_{ls} = \Ab^{\dagger}\yb$.

Computing $\Ab^{\dagger}$ to solve \eqref{OLS} directly is intractable for large $m$, as it requires computing the inverse of $\Ab^T\Ab$. Instead, we use gradient methods to get \textit{approximate} solutions, e.g. SD or CG, which require less operations, and can be done distributively. One could of course use second-order methods, e.g. such as Newton–Raphson, Gauss-Newton, Quasi-Newton, BFGS.

We briefly review steepest descent. When considering a minimization problem with a convex differentiable objective function $\phi\colon\Theta\to\R$ over an open constrained set $\Theta\subseteq \R^m$, {as in \eqref{OLS}}, we select an initial $\theta^{(0)}\in\Theta$ and repeat:
$$ \theta^{(t+1)}=\theta^{(t)}-\xi_t\cdot\nabla_{\theta}\phi(\theta^{(t)}), \quad \text{ for } t=1,2,3,...$$
until a termination criterion is met. The criterion may depend on the problem we are looking at, and the parameter $\xi_t$ is the step-size, which may be adaptive or fixed. In our experiments, we use backtracking line search to determine $\xi_t$.

\section{Approximation Algorithms}
\label{Appr_alg_sec}

\subsection{Proposed Inverse Algorithm}
\label{pr_inv_alg_subsec}

$\ind$ Our goal is to estimate $\Ab^{-1}=\big[\bb_1 \ \cdots \ \bb_n \big]$, for $\Ab$ a square matrix of order $n$. A key property to note is
$$ \Ab\Ab^{-1} = \Ab\big[\bb_1 \ \cdots \ \bb_n \big] = \big[\Ab\bb_1 \ \cdots \ \Ab\bb_n \big] = \bold{I}_n $$
which implies that $\Ab\bb_i=\eb_i$ for all $i\in\N_n\coloneqq\{1,\cdots,n\}$, where $\eb_i$ are the standard basis column vectors. Assume for now that we use any black-box least squares solver to estimate
\begin{equation}
\label{inv_LS}
  \hat{\bb}_i = \arg\min_{\bb\in\R^n} \Big{\{f_i(\bb)\coloneqq\|\Ab\bb-\eb_i\|_2^2}\Big\}
\end{equation}
which we call $n$ times, to estimate $\widehat{\Ab^{-1}} = \big[ \hat{\bb}_1 \ \cdots \ \hat{\bb}_n \big]$. This approach may be viewed as solving
\begin{equation*}
\label{inv_LS_F}
  \widehat{\Ab^{-1}} = \arg\min_{\ \ \Bb\in\R^{n\times n}}\left\{\|\Ab\Bb-\bold{I}_n\|_F^2\right\}.
\end{equation*}
Alternatively, one could estimate the rows of $\Ab^{-1}$. Algorithm \ref{inv_alg} shows how this can be performed by a single server.

\begin{algorithm}[h]
\label{inv_alg}
\SetAlgoLined
\KwIn{nonsingular $\Ab\in\R^{n\times n}$}
\KwOut{estimate $\widehat{\Ab^{-1}}$ of $\Ab$'s inverse}
  \For{i=1 to n}
  {
    solve $\hat{\bb}_i = \arg\min_{\bb\in\R^n} \left\{\|\Ab\bb-\eb_i\|_2^2\right\}$ 
  }
 \Return $\widehat{\Ab^{-1}} \gets \big[ \hat{\bb}_1 \ \cdots \ \hat{\bb}_n \big]$
 \caption{Estimating $\Ab^{-1}$}
\end{algorithm}

In the case where SD is used to estimate each column, the overall operation count is $\ow(nTn^2)$, where $T$ is the average number of iterations used per column estimation. The number of iterations can be determined by the termination criterion, e.g. the criterion $f_i(\bbh^{(t)})-f_i(\bb^{\star})\leq\epsilon$ is guaranteed to be satisfied after $T=\ow(\log(1/\epsilon))$ iterations \cite{BV04}. The overall error of our estimate may be quantified as
\begin{itemize}
  \item $\text{err}_{\ell_2}(\widehat{\Ab^{-1}}) \coloneqq \|\widehat{\Ab^{-1}}-\Ab^{-1}\|_2^2$
  \item $\text{err}_F(\widehat{\Ab^{-1}}) \coloneqq \|\widehat{\Ab^{-1}}-\Ab^{-1}\|_F^2$
  \item $\text{err}_{\text{r}F}(\widehat{\Ab^{-1}}) \coloneqq \frac{\|\widehat{\Ab^{-1}}-\Ab^{-1}\|_F^2}{\|\Ab^{-1}\|_F^2} = \frac{\sum\limits_{i=1}^{n} \|\Ab\hat{\bb}_i-\eb_i\|_2^2}{\|\Ab^{-1}\|_F^2}$
\end{itemize}
which we refer to as the \textit{$\ell_2$-error}, \textit{Frobenius-error} and \textit{relative Frobenius-error} respectively. The corresponding pseudoinverse approximation errors are defined accordingly.

If there are $n$ servers then Algorithm \ref{inv_alg} can be used to estimate $\Ab^{-1}$ in a distributive manner, where the workers each compute a $\bbh_i$ in parallel. The expected runtime assuming no delays in this case is therefore $\ow(T_{\text{max}}n^2)$, for $T_{\text{max}}$ the maximum number of iterations of SD used by the workers to estimate their respective column. This is depicted in Figure 1.

\begin{figure}[h]
  \centering
  \label{inv_appr_recovery}
    \includegraphics[scale=.18]{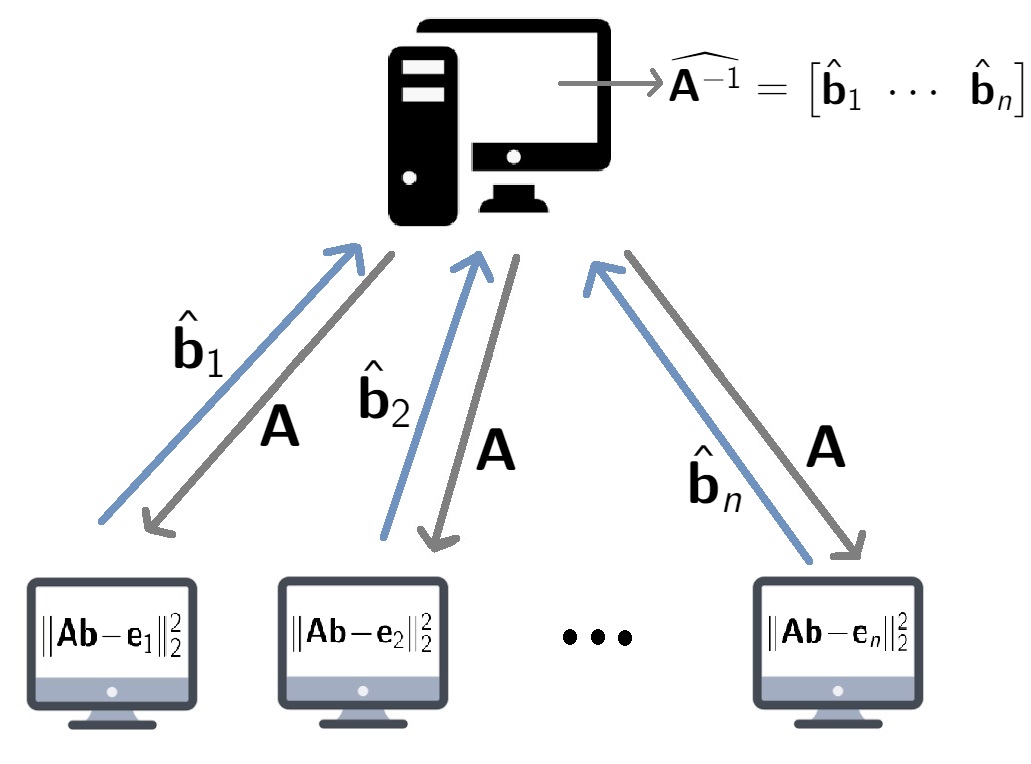}
    \caption{Communication schematic, with the recovery of $\widehat{\Ab^{-1}}$.}
\end{figure}

In contrast to the algorithms in \cite{MHGB18,LLA16,XMA14,QQSV01,QQSV98,YLX13,LKV96,BF88}, Algorithm \ref{inv_alg} when carried out distributively, requires the central and worker servers to communicate only once, after the workers are provided with $\Ab$. Furthermore, the communication load is optimal, in the sense that the workers send the same amount of information once, and each send the minimal amount of matrix entries under the constraint which is determined by the number of stragglers the scheme tolerates. Workers could also be assigned to estimate multiple columns. This is discussed further in Section \ref{Str_pr_sec}.

\subsection{Frobenius and Relative Frobenius Errors}
\label{rel_err_subsec}

In order to bound $\text{err}_{\text{r}F}(\widehat{\Ab^{-1}})=\frac{\|\widehat{\Ab^{-1}}-\Ab^{-1}\|_F^2}{\|\Ab^{-1}\|_F^2}$, we first upper bound the numerator and then lower bound the denominator. Since $\|\Ab^{-1}-\widehat{\Ab^{-1}}\|_F^2=\sum_{i=1}^n\|\Ab^{-1}\eb_i-\hat{\bb}_i\|_2^2$, bounding the numerator reduces to bounding $\|\Ab^{-1}\eb_i-\hat{\bb}_i\|_2^2$ for all $i\in\N_n$. This calculation is straightforward
\begin{align}
\label{upp_bd_numer}
  \|\Ab^{-1}\eb_i-\hat{\bb}_i\|_2^2 &\overset{\Diamond}{\leqslant} 2\|\Ab^{-1}\eb_i\|_2^2+2\|\hat{\bb}_i\|_2^2\notag\\
  &\overset{\$}{\leqslant} 2\|\Ab^{-1}\|_2^2\cdot\|\eb_i\|_2^2+2\|\hat{\bb}_i\|_2^2\notag\\
  &= 2/\sigma_{\text{min}}(\Ab)^2+2\|\hat{\bb}_i\|_2^2
\end{align}
where in $\Diamond$ we use the fact that $\|\ub+\vb\|_2^2\leqslant2\|\ub\|_2^2+2\|\vb\|_2^2$, and in $\$$ the submultiplicativity of the $\ell_2$-norm. For the denominator, by the definition of the Frobenius norm
\begin{equation}
\label{low_bd_denom}
  \|\Ab^{-1}\|_F^2=\sum_{i=1}^n\frac{1}{\sigma_i(\Ab)^2} \geq \frac{n}{\sigma_{\text{max}}(\Ab)^2}.
\end{equation}
By combining \eqref{upp_bd_numer} and \eqref{low_bd_denom} we get
\begin{align*}
  \text{err}_{\text{r}F}(\widehat{\Ab^{-1}}) &\leq 2\cdot\frac{n/\sigma_{\text{min}}(\Ab)^2+\sum_{i=1}^n\|\hat{\bb}_i\|_2^2}{n/\sigma_{\text{max}}(\Ab)^2}\\
  &= 2\kappa_2(\Ab)^2+\frac{2\sigma_{\text{max}}(\Ab)^2}{n}\cdot\sum_{i=1}^n\|\hat{\bb}_i\|_2^2
\end{align*}
for $\kappa_2(\Ab)=\|\Ab\|_2\|\Ab^{-1}\|_2$ the condition number of $\Ab$. This shows a dependency on the estimates, and is an additive error bound in terms of the problem's condition number.

Proposition \ref{rel_err_bound_prop} gives an error bound when using SD to solve \eqref{inv_LS}, with termination criterion $\|\nabla f_i(\bb^{(k)})\|_2\leq\epsilon$.

\begin{Prop}
\label{rel_err_bound_prop}
  For $\Ab\in\R^{n\times n}$ full-rank, we have $\mathrm{err}_{F}(\widehat{\Ab^{-1}})\leq \frac{n\epsilon^2/2}{\sigma_{\mathrm{min}}(\Ab)^4}$ and $\mathrm{err}_{\mathrm{r}F}(\widehat{\Ab^{-1}})\leq \frac{n\epsilon^2/2}{\sigma_{\mathrm{min}}(\Ab)^2}$, when using SD to solve \eqref{inv_LS} with termination criteria $\|\nabla f_i(\bb^{(k)})\|_2\leq\epsilon$ for all $i\in\N_n$.
\end{Prop}

\begin{proof}
Recall that for a strongly-convex function with strong-convexity parameter $\mu$, we have the following optimization gap \cite[Section 9.1.2]{BV04}
\begin{equation}
\label{opt_gap}
  f_i(\bb)-f_i(\bb_{ls}^{\star})\leq\frac{1}{2\mu}\cdot\|\nabla f_i(\bb)\|_2 \ .
\end{equation}
For $\Ab$ full-rank in \eqref{inv_LS}, the constant is $\mu=\sigma_{\text{min}}(\Ab)^2$. By fixing $\epsilon=\sqrt{2\sigma_{\text{min}}(\Ab)^2\alpha}$, we have $\alpha=\frac{1}{2}\cdot\left(\frac{\epsilon}{\sigma_{\text{min}}(\Ab)}\right)^2$. Thus by \eqref{opt_gap} and our termination criterion
$$ \|\nabla f_i(\bb)\|_2\leq\sqrt{2\sigma_{\text{min}}(\Ab)^2\alpha} \quad \implies \quad f_i(\bb)-f_i(\bb_{ls}^{\star})\leq\alpha\ , $$
so when solving \eqref{inv_LS} we get
$$ f_i(\bb)-f_i(\bb_{ls}^{\star}) = f_i(\bb)-0 = \|\Ab\hat{\bb}_i-\eb_i\|_2^2 $$
hence
\begin{equation}
\label{summand_obj_bound}
  \|\Ab\hat{\bb}_i-\eb_i\|_2^2 \leq \frac{1}{2}\cdot\left(\frac{\epsilon}{\sigma_{\text{min}}(\Ab)}\right)^2
\end{equation}
for all $i\in\N_n$. We want an upper bound for each summand $\|\Ab^{-1}\eb_i-\hat{\bb}_i\|_2^2$ of the numerator of $\text{err}_{\text{r}F}(\widehat{\Ab^{-1}})$
\begin{align}
\label{deriv_rF_bd}
  \|\Ab^{-1}\eb_i-\hat{\bb}_i\|_2^2 &= \|\Ab^{-1}(\eb_i-\Ab\hat{\bb}_i)\|_2^2\notag\\
  &\leq \|\Ab^{-1}\|_2^2\cdot\|\eb_i-\Ab\hat{\bb}_i\|_2^2\notag\\
  &\overset{\sharp}{\leq} \|\Ab^{-1}\|_2^2\cdot\frac{1}{2} \cdot \left(\frac{\epsilon}{\sigma_{\text{min}}(\Ab)}\right)^2\notag\\
  &= \epsilon^2/\left(2\sigma_{\text{min}}(\Ab)^4\right)
\end{align}
where $\sharp$ follows from \eqref{summand_obj_bound}, thus $\text{err}_{F}(\widehat{\Ab^{-1}})\leq\frac{n\epsilon^2}{2}$. Substituting this into the definition of $\text{err}_{\text{r}F}(\widehat{\Ab^{-1}})$ gives us
$$ \text{err}_{\text{r}F}(\widehat{\Ab^{-1}}) \leq \frac{\|\Ab^{-1}\|_2^2}{\|\Ab^{-1}\|_F^2} \cdot\frac{n}{2} \cdot \left(\frac{\epsilon}{\sigma_{\text{min}}(\Ab)}\right)^2 \overset{\ddagger}{\leq}\frac{n\epsilon^2/2}{\sigma_{\mathrm{min}}(\Ab)^2} $$
where $\ddagger$ follows from the fact that $\|\Ab^{-1}\|_2^2\leq\|\Ab^{-1}\|_F^2$.
\end{proof}

In the experiments of Section \ref{exper_sec}, we verify that Proposition \ref{rel_err_bound_prop} holds for Gaussian random matrices. The dependence on $1/\sigma_{\text{min}}(\Ab)$ is an artifact of using gradient methods and the underlying problems \eqref{inv_LS}, since the error will be multiplied by $\|\Ab^{-1}\|_2^2$. In theory, this can be annihilated if one runs the algorithm on $d\Ab$ for $d\approx1/\sigma_{\text{min}}(\Ab)$, and then multiply the final result by $d$. This is a way of preconditioning SD. In practice, the scalar $d$ should not be selected much larger than $1/\sigma_{\text{min}}(\Ab)$, as it would result in $\widehat{\Ab^{-1}}\approx\bold{0}_{n\times n}$.

\subsection{Proposed Pseudoinverse Algorithm}
\label{pr_pinv_alg_subsec}

Just like the inverse, the pseudoinverse of a matrix also appears in a variety of applications and computations in numerous fields. Computing the pseudoinverse is even more cumbersome, as itself requires computing an inverse. For this subsection, we consider $\Ab\in\R^{n\times m}$ for $n>m$, of rank $m$.

One could naively attempt to modify Algorithm \ref{inv_alg} in order to retrieve $\Ab^{\dagger}$ such that $\Ab^{\dagger}\Ab=\Ib_m$, by approximating the rows of $\Ab^{\dagger}$. This would \textit{not} work, as the underlying optimization problems would not be strictly convex. Instead, we use Algorithm \ref{inv_alg} to estimate $\Bb^{-1}\coloneqq(\Ab^T\Ab)^{-1}$, and then multiply the estimate $\widehat{\Bb^{-1}}$ by $\Ab^T$. The multiplication may be done by the workers (who will estimate the rows of $\Ab^{\dagger}$) or the central server. For coherence, we present the case where the workers carry out the multiplication.

The additional operation count compared to Algorithm \ref{inv_alg}, is $\ow(m^2n)$ for constructing $\Bb$, and $\ow(mn)$ per worker for multiplying their approximation $\hat{\cb}_i$ by $\Ab^T$.

\begin{algorithm}[h]
\label{pinv_alg}
\SetAlgoLined
\KwIn{full-rank $\Ab\in\R^{n\times m}$ where $n>m$}
\KwOut{estimate $\widehat{\Ab^{\dagger}}$ of $\Ab$'s pseudoinverse}
  $\Bb\gets \Ab^T\Ab$\\
  \For{i=1 to m}
  {
    solve $\hat{\cb}_i = \arg\min_{\cb\in\R^{1\times m}} \Big\{g_i(\cb)\coloneqq\|\cb\Bb-\eb_i^T\|_2^2\Big\}$\\ 
    $\hat{\bb}_i\gets \hat{\cb}_i\cdot\Ab^T$
  }
 \Return $\widehat{\Ab^{\dagger}} \gets \left[ \hat{\bb}_1^T \ \cdots \ \hat{\bb}_m^T \right]^T$
 \caption{Estimating $\Ab^{\dagger}$}
\end{algorithm}

\begin{Cor}
\label{rel_err_bound_cor}
  For full-rank $\Ab\in\R^{n\times m}$ with $n>m$, we have $\mathrm{err}_{F}(\widehat{\Ab^{\dagger}})\leq m\cdot\left(\frac{\epsilon\cdot\kappa_2(\Ab)}{\sqrt{2}\sigma_{\mathrm{min}}(\Ab)}\right)^2$ and $\mathrm{err}_{\mathrm{r}F}(\widehat{\Ab^{\dagger}})\leq \frac{m\epsilon^2\cdot\kappa_2(\Ab)^2}{2}$ when using SD to solve the subroutine optimization problems of Algorithm \ref{pinv_alg}, with termination criteria $\|\nabla g_i(\cb^{(k)})\|_2\leq\epsilon$.
\end{Cor}

\begin{proof}
Similarly to the derivation of \eqref{deriv_rF_bd}, we get
$$ \|\Bb^{-1}\eb_i-\hat{\cb}_i^T\|_2^2 \leq \left(\frac{\epsilon/\sqrt{2}}{\sigma_{\mathrm{min}}(\Bb)}\right)^2 = \left(\frac{\epsilon/\sqrt{2}}{\sigma_{\mathrm{min}}(\Ab)^2}\right)^2 \eqqcolon \delta \ . $$
Denote by $\Ab_{(i)}^{\dagger}$ the $i^{th}$ row of  $\Ab^{\dagger}$. The above bound implies that for each summand of the Frobenius error;
$\|\hat{\bb}_i-\Ab^{\dagger}_{(i)}\|_2^2=\|\hat{\cb}_i\Ab^T-\eb_i^T\cdot\Bb^{-1}\Ab^T\|_2^2$, we have $\|\hat{\bb}_i-\Ab^{\dagger}_{(i)}\|_2^2\leq\delta\|\Ab^T\|_2^2$. Summing the right hand side $m$ times and using the fact that $1/\sigma_{\text{min}}(\Ab)^2=\|\Ab^{\dagger}\|_2^2\leq\|\Ab^{\dagger}\|_F^2$, completes the proof.

\end{proof}

\section{Robust to Stragglers}
\label{Str_pr_sec}  

We now discuss coded computing, and give a binary fractional repetition scheme \cite{ERK10,TLDK17,CMH20} which makes Algorithms \ref{inv_alg} and \ref{pinv_alg} robust to stragglers. We present the proposed scheme for Algorithm \ref{inv_alg}. While there is extensive literature on matrix-matrix, matrix-vector multiplication \cite{LLPPR18,YMAA17,LSR17,YMAA20,YA20,FJHDCG17,DFHJCG19,FC19,SHN19,CPH20c}, and computing the gradient \cite{TLDK17,HASH17,RTTD17,OGU19,CMH20,YA18,CP18,CPE17,WCP19,BWE19,WLS19,KKR19,HYKM19,CPH20a} in the presence of stragglers, there is limited work on computing or approximating the inverse of a matrix \cite{YGK17}.

For our algorithms, any coded computing scheme in which the workers compute an encoding of partitions of the resulting computation $\Eb=\big[E_1 \ \cdots \ E_k \big]$ (i.e. matrix product) could be utilized. It is crucial that the encoding takes place on the computed tasks $\{E_i\}_{i=1}^k$ in the scheme, and not the assigned data or partitions of the matrices that are being computed over, otherwise Algorithms \ref{inv_alg} and \ref{pinv_alg} would not return the correct results. This corresponds to a partitioning $\widehat{\Ab^{-1}} = \big[\Acal_1 \ \cdots \ \Acal_k \big]$ of our approximation. We resort to \textit{gradient coding} (GC) rather than \textit{coded matrix multiplication}, as most multiplication schemes encode the matrices a priori.

Assume that each computational task $\{\Acal_i\}_{i=1}^k$ is comprised of $\texT=\frac{n}{k}$ distinct but consecutive approximations of \eqref{inv_LS}, i.e.
$$ \Acal_i = \big[\bbh_{(i-1)\texT+1}\ \cdots \ \bbh_{i\texT}\big]\in\R^{n\times\texT} \quad \text{ for each } i=1,\cdots,k\ . $$
If $k\nmid n$, we allocate $\tilde{\texT}=\floor{\frac{n}{k}}$ many approximate column vectors $\bbh_j$ to $\tilde{k}=\text{rem}(n,k)$ of the partitions $\Acal_i$, and $\tilde{\texT}+1$ approximate vectors to the remaining $k-\tilde{k}$ partitions.

To simplify our presentation, we assume that $(s+1)\mid N$ and $k\mid n$ for $s$ the maximum number of stragglers, and let $\ell=\frac{N}{s+1}$ to be equal to $k$ (as was assumed in \cite{TLDK17}). Moreover, we assume the workers are \textit{homogeneous} (i.e. they have the same computational power), and therefore equal computational loads are assigned to them. The analysis and modification to the encoding algorithm for heterogeneous workers in \cite{CMH20} also applies to the coded scheme we present below.

We adapt the GC scheme from \cite{CMH20} which relies on the pigeonhole principle, to our setting, when Algorithms \ref{inv_alg} and \ref{pinv_alg} are deployed in a distributed manner. The correspondence between (i) the GC scheme and (ii) distributed matrix inversion, is that the partial gradient $\gb_i$ in (i) corresponds to $\Acal_i$ in (ii). Considering $N$ workers, by requesting $s+1$ workers to compute $\Acal_i$ for each $i\in\N_k$, as long as $f\coloneqq N-s$ workers respond, the approximation scheme tolerates $s$ stragglers and $\widehat{\Ab^{-1}}$ is recoverable. We allocate to the workers with consecutive indices the same computation, i.e. the workers with indices $\{(i-1)\cdot(s+1),\cdots,i\cdot(s+1)\}$ each compute $\Acal_i$, for each $i$. In such a scenario, the communication load per worker will scale according to the number of columns which comprise each $\Acal_i$, when compared to the GC scheme. It is crucial that there is no overlap between the assigned columns of workers with congruent indices $\bmod (s+1)$\footnote{Recall that $a$ is congruent to $b$ if $a\equiv b\bmod(s+1)$.}.

Once $f$ workers respond, we are guaranteed by the pigeonhole that at least one complete residue system of indices has responded. That is, there is at least one integer $l\in\N_{s+1}$ for which all workers with index congruent to $l\bmod(s+2)$ have responded \cite{CMH20}. The central server knows the index set $\I$ of the $f$ responsive workers, and can therefore determine an integer $l$ which corresponds to a complete residue system\footnote{More than one may occur, and we are guaranteed that at least one always occurs, as long as no less than $N-s$ workers respond.}. Note that we are actually working $\bmod(s+1)$, but since we indexed the workers starting from 1, we describe everything $\bmod(s+2)$ and do not consider the congruence class which is equivalent to $0\bmod(s+1)$.

We set up a framework analogous to that proposed for GC in \cite{TLDK17}, where the objective is to design an encoding matrix $\Gb$ and a decoding vector $\ab_{\I}$ which is determined by the responsive workers' index set $\I$, so that $\ab_{\I}^T\Gb=\bold{1}_{1\times k}$ regardless of $\I$. The cardinality of $\I$ is always $f$. The objective now is to concatenate the computations $\Acal_i$ instead of summing them (in GC the objective is to recover $\gb=\sum_{i=1}^k\gb_i$). Therefore, we want to construct matrices $\Gbt\in\{0,1\}^{n(s+1)\times(s+1)}$ and $\abt_{\I}\in\{0,1\}^{n(s+1)\times(s+1)}$ that satisfy $\abt_{\I}^T\Gbt=\Ib_n$, for \textit{any} of the ${{N}\choose{f}}$ possible index sets $\I$. The encoding step of our scheme corresponds to
\begin{equation}
\label{enc_identity}
  \overbrace{\big(\Ib_k\otimes\bold{1}_{(s+1)\times 1}\otimes\Ib_{\texT}\big)}^{\Gbt\in\{0,1\}^{k\texT(s+1)\times k\texT}} \cdot \overbrace{\begin{bmatrix} \Acal_1^T \\ \vdots \\ \Acal_k^T \end{bmatrix}}^{(\widehat{\Ab^{-1}})^T} = \begin{bmatrix} \Acal_1^T \\ \vdots \\ \Acal_1^T \\ \vdots \\ \Acal_k^T \\ \vdots \\ \Acal_k^T \end{bmatrix}
\end{equation}
where the transpose of each submatrix $\Acal_i$ appears $s+1$ times along the rows of the encoding $\Gbt(\widehat{\Ab^{-1}})^T\in\R^{(s+1)n\times n}$. Each block $\Acal_i$ corresponds to one of the $s+1$ workers which is asked to compute $\Acal_i$. Hence, there are $k(s+1)$ blocks in the encoding $\Gbt(\widehat{\Ab^{-1}})^T$ \eqref{enc_identity}. The encoding matrix $\Gbt$ reveals the task allocation which is applied to the computations $\{\Acal_i\}_{i=1}^k$ by $\Gbt(\widehat{\Ab^{-1}})^T\in\R^{(s+1)n\times n}$. Algorithmically, $\Gbt$ is constructed as described in Algorithm \ref{enc_Gt}.

\begin{algorithm}[h]
\label{enc_Gt}  
\SetAlgoLined
  \KwIn{design parameters $N$, $s$ and $\texT=n/k$}
  \KwOut{$\Gbt\in\{0,1\}^{n(s+1)\times n}$}
  \textbf{Initialize:} {\small$\ell=\frac{N}{s+1}$ and $\Gbt\gets\{0,1\}^{N\times \ell}$ \Comment{assume $k=\ell$}}\\
  \For{$j=1$ to $\ell$}
    {$\Gbt[(j-1)\cdot(s+1)+1:j\cdot(s+1),j]=\bold{1}_{(s+1)\times1}$}
  \Return $\Gbt\gets\Gbt\otimes\Ib_{\texT}$
\caption{Encoding matrix $\Gbt$}
\end{algorithm}

The decoding matrix $\abt_{\I}\in\{0,1\}^{\frac{Nn}{k}\times n}$ constructed by Algorithm \ref{determine_aI_matrix}, is determined by the integer $l$ which corresponds to a complete residue system in $\I$.

\begin{algorithm}[h]
\label{determine_aI_matrix}
\SetAlgoLined
  \KwIn{index set $\I$ of  $f$ fastest workers, and design parameters $n$, $s$ and $\texT=n/k$}
  \KwOut{$\abt_{\I}\in\{0,1\}^{n(s+1)\times n}$}
  \textbf{Initialize:} {\small$\ell=\frac{N}{(s+1)}$ and $\abt\gets\{0,1\}^{N\times \ell}$ \Comment{assume $k=\ell$}}\\
  Determine the integer $l$ which corresponds to a complete residue system in $\I$\\
  \For{$j=1$ to $\ell$}
    {$\abt[(j-1)\cdot(s+1)+l,j]=1$}
  \Return $\abt_{\I}\gets\abt\otimes\Ib_{\texT}$
\caption{Decoding matrix $\abt_{\I}$}
\end{algorithm}

\begin{Prop}
The coding scheme comprised of $\Gbt$ and $\abt_\I$ based on Algorithms \ref{enc_Gt} and \ref{determine_aI_matrix} is robust to $s$ stragglers.
\end{Prop}

\begin{proof}
We are guaranteed that always one complete residue system is present in $\I$, as long as $f=N-s$ workers respond. Hence, $\abt_\I$ can always be constructed by Algorithm \ref{determine_aI_matrix}. 

From the encoding of the computations by $\Gbt$ \eqref{enc_identity}, we receive $f$ blocks $\Acal_j$. There may be repetitions, but at least one of each $\{\Acal_i\}_{i=1}^k$ is received. The decoding matrix $\abt_\I$ only retrieves those computations corresponding to the workers with indices congruent to $l\bmod(s+1)$, and ignores the rest. Therefore $\abt_\I^T\Gbt=\Ib_n$ for \textit{any} $\I$, and $\abt_\I^T\cdot\big(\Gbt\cdot(\widehat{\Ab^{-1}})^T\big)=(\widehat{\Ab^{-1}})^T$.
\end{proof}

We point out that this scheme is numerically stable, as the encoding and decoding steps do not introduce instability when they are applied. The scheme is also optimal according to the analysis in \cite{CMH20}, where it is assumed that the encoding takes place on the computed tasks and \textit{not} the data. Analogously, we encode $\{\Acal_i\}_{i=1}^k$ and \textit{not} $\Ab$. Furthermore, since the encoding can be viewed as a task assignment, the computations can be sent in a sequential manner, and then terminate once a complete residue system of worker indices is received. A similar idea appears in \cite{OGU19}, in a different context.


We give a simple example of our scheme, to help visualize the encoding task assignments and the decoding. Let $N=8$, $s=1$, thus $\ell=k=4$, and $n$ be arbitrary, with 
$\texT=n/k$. By $\bold{0}_{\texT}$ and $\Ib_\texT$ we denote the zero square and identity matrices respectively, of size $\texT\times\texT$.  Each of the $s+1=2$ congruence classes are represented by a different color and font. For $l=2$, i.e. $\{2,4,6,8\}\subseteq\I$, the encoding decoding pair is
$$ \Gbt = \begin{bmatrix} {\cyan \Ib_{\texT}} & & & \\ {\purple \mathbb{I}_{\mathbb{T}}} & & & \\ & {\cyan \Ib_{\texT}} & &  \\ & {\purple \mathbb{I}_{\mathbb{T}}} & & \\ & & {\cyan \Ib_{\texT}} & \\ & & {\purple \mathbb{I}_{\mathbb{T}}} & \\ & & & {\cyan \Ib_{\texT}} \\ & & & {\purple \mathbb{I}_{\mathbb{T}}} \end{bmatrix} \ , \ \abt_{\I} = \begin{bmatrix} {\cyan \bold{0}_{\texT}} & & & \\ {\purple \mathbb{I}_{\mathbb{T}}} & & & \\ & {\cyan \bold{0}_{\texT}} & &  \\ & {\purple \mathbb{I}_{\mathbb{T}}} & & \\ & & {\cyan \bold{0}_{\texT}} & \\ & & {\purple \mathbb{I}_{\mathbb{T}}} & \\ & & & {\cyan \bold{0}_{\texT}} \\ & & & {\purple \mathbb{I}_{\mathbb{T}}} \end{bmatrix}
$$
where both matrices are of the same size. From this example, it is also clear that $\abt_\I$ is in fact the restriction of $\Gbt$ to the workers corresponding to the $l^{th}$ congruence class, when $(s+1)|N$.

The same scheme can be applied when estimating the pseudoinverse $\Ab^{\dagger}$. The only difference is that the computations $\widehat{\Ab^{\dagger}} = \big[\Acal_1 \ \cdots \ \Acal_k \big]$ according to Algorithm \ref{pinv_alg} are now allocated to the workers, hence the encoding is $\Gbt\widehat{\Ab^{\dagger}}$.


\section{Experiments}
\label{exper_sec}

The accuracy of the proposed algorithms was tested on randomly generated matrices, using both SD and CG \cite{TB97} for the subroutine optimization problems. The depicted results are averages of 20 runs, with termination criteria $\|\nabla f_i(\bb^{(k)})\|_2\leq \epsilon$ for SD and $\|\bb_i^{(k)}-\bb_i^{(k-1)}\|_2\leq \epsilon$ for CG, for the given $\epsilon$ accuracy parameters. The criteria for $\widehat{\Ab^{\dagger}}$ were analogous. We considered $\Ab\in \R^{100\times 100}$ and $\Ab\in\R^{100\times 50}$. The error subscripts represent $\mathscr{A}=\{\ell_2,{F},\text{r}F\}$, $\mathscr{N}=\{\ell_2,F\}$, $\mathscr{F}=\{F,\text{r}F\}$. We note that significantly fewer iterations took place when CG was used for the same $\epsilon$. Thus, there is a trade-off between accuracy and speed when using SD vs. CG.

\begin{center}
\begin{tabular}{ |p{.5cm}||p{1.05cm}|p{1.1cm}|p{1.1cm}|p{1.1cm}|p{1.1cm}| }
\hline
\multicolumn{6}{|c|}{Average $\widehat{\Ab^{-1}}$ errors for $\Ab\sim50\cdot \mathcal{N}(0,1)$ --- SD} \\
\hline
$\epsilon$ & $10^{-1}$ & $10^{-2}$ & $10^{-3}$ & $10^{-4}$ & $10^{-5}$ \\
\hline
$\text{err}_{\mathscr{A}}$ & {\small$\ow(10^{-2})$} & {\small$\ow(10^{-5})$} & {\small$\ow(10^{-7})$} & {\small$\ow(10^{-9})$} & {\small$\ow(10^{-12})$} \\
\hline
\end{tabular}
\end{center}

\begin{center}
\begin{tabular}{ |p{.5cm}||p{1.05cm}|p{1.05cm}|p{1.05cm}|p{1.1cm}|p{1.1cm}| }
\hline
\multicolumn{6}{|c|}{{\small Mean $\widehat{\Ab^{-1}}$ errors, $\Ab=\bold{M}+\bold{M}^T$, $\bold{M}\sim25\cdot \mathcal{N}(0,1)$ --- CG}} \\
\hline
$\epsilon$ & $10^{-3}$ & $10^{-4}$ & $10^{-5}$ & $10^{-6}$ & $10^{-7}$ \\
\hline
$\text{err}_{\mathscr{N}}$ & {\small$\ow(10^{-3})$} & {\small$\ow(10^{-5})$} & {\small$\ow(10^{-8})$} & {\small$\ow(10^{-11})$} & {\small$\ow(10^{-12})$} \\
$\text{err}_{\text{r}F}$ & {\small$\ow(10^{-3})$} & {\small$\ow(10^{-5})$} & {\small$\ow(10^{-7})$} & {\small$\ow(10^{-10})$} & {\small$\ow(10^{-12})$} \\
\hline
\end{tabular}
\end{center}

\begin{center}
\begin{tabular}{ |p{.5cm}||p{1.05cm}|p{1.05cm}|p{1.05cm}|p{1.15cm}|p{1.15cm}| }
\hline
\multicolumn{6}{|c|}{Average $\widehat{\Ab^{\dagger}}$ errors for $\Ab\sim \mathcal{N}(0,1)$ --- SD} \\
\hline
$\epsilon$ & $10^{-1}$ & $10^{-2}$ & $10^{-3}$ & $10^{-4}$ & $10^{-5}$ \\
\hline
$\text{err}_{\ell_2}$ & {\small$\ow(10^{-4})$} & {\small$\ow(10^{-6})$} & {\small$\ow(10^{-8})$} & {\small$\ow(10^{-10})$} & {\small$\ow(10^{-12})$} \\
$\text{err}_{\mathscr{F}}$ & {\small$\ow(10^{-5})$} & {\small$\ow(10^{-7})$} & {\small$\ow(10^{-9})$} & {\small$\ow(10^{-11})$} & {\small$\ow(10^{-13})$} \\
\hline
\end{tabular}
\end{center}

\begin{center}
\begin{tabular}{ |p{.5cm}||p{1.05cm}|p{1.05cm}|p{1.05cm}|p{1.15cm}|p{1.15cm}| }
\hline
\multicolumn{6}{|c|}{Average $\widehat{\Ab^{\dagger}}$ errors for $\Ab\sim \mathcal{N}(0,1)$ --- CG} \\
\hline
$\epsilon$ & $10^{-3}$ & $10^{-4}$ & $10^{-5}$ & $10^{-6}$ & $10^{-7}$ \\
\hline
$\text{err}_{\ell_2}$ & {\small$\ow(10^{-4})$} & {\small$\ow(10^{-6})$} & {\small$\ow(10^{-8})$} & {\small$\ow(10^{-10})$} & {\small$\ow(10^{-12})$} \\
$\text{err}_{\mathscr{F}}$ & {\small$\ow(10^{-2})$} & {\small$\ow(10^{-3})$} & {\small$\ow(10^{-8})$} & {\small$\ow(10^{-10})$} & {\small$\ow(10^{-12})$} \\
\hline
\end{tabular}
\end{center}


\bibliographystyle{IEEEtran}
\bibliography{refs}

\end{document}